\documentclass{article}
\usepackage{graphicx}
\usepackage{amsmath}
\usepackage{amsfonts}
\usepackage{amssymb}
\newtheorem{theorem}{Theorem}

\newtheorem{corollary}[theorem]{Corollary}

\newenvironment{proof}[1][Proof]{\textbf{#1.} }{\ \rule{0.5em}{0.5em}}

\begin{document}

\title{Groupoid Approach to Quantum Projective Spaces}
\author{Albert Jeu-Liang Sheu\thanks{Partially supported by NSF Grant DMS-9623008.}\\Department of Mathematics\\University of Kansas\\Lawrence, KS 66045\\U. S. A.}
\date{}
\maketitle

\section{Introduction}

The theory of groupoid C*-algebras \cite{Re} has proved to be a very useful
tool and provide a nearly universal context for the study of operator
algebras, especially in connection with geometric structures. There have been
many successful examples of such uses of groupoid C*-algebras
\cite{Co,CuM,MRe,SSU,Sh:rd}.

Recently interesting and important examples of C*-algebras of quantum groups
$G_{q}$ and quantum spaces $M_{q}$ arise from the quantization
\cite{D,FRT,So,Wo:ts,Wo:cm,Po,Ri:cq,Ri:dq,Ri:dqar,VaSo:af,Sh:qp,N:dq} of
Poisson Lie groups $G$ and their Poisson homogeneous spaces $M=H\backslash G$
with $H $ a Poisson Lie subgroup of $G$. Known results \cite{So,Sh:qp} show
that they are closely related to the underlying singular symplectic foliation
\cite{We,LuWe:plg}. It is shown in \cite{Sh:cqg} that the C*-algebra of
compact quantum groups can be effectively described and studied in the context
of groupoid C*algebras with the groupoid summarizing the underlying singular
symplectic foliation. In fact, the structures of the algebras of quantum
groups $SU\left(  n\right)  _{q}$ and quantum spheres $\mathbb{S}_{q}^{2n-1}$
are analyzed in detail in \cite{Sh:cqg,Sh:qs}.

More recently, nonstandard quantum complex projective spaces $\mathbb{CP}%
_{q,c}^{n}$ have been studied by Dijkhuizen and Noumi \cite{DiNo}. In this
paper, we show how the C*-algebras of quantum complex projective spaces
(standard or nonstandard) are related to groupoids.

\section{Quantum groups $SU\left(  n\right)  _{q}$}

In this section, we briefly recall the notations and results needed from
\cite{Sh:cqg}. For more complete information, we refer to the paper
\cite{Sh:cqg}.

First we recall that the transformation group groupoid ${\mathbb{Z}}^{m}%
\times\overline{{\ \mathbb{Z}}}^{m}$ (with ${\mathbb{Z}}^{m}$ acting on
$\overline{{\mathbb{Z}}}^{m}$ by translation) when restricted to the positive
cone $\overline{{\mathbb{Z}}}_{\geq}^{m}$ gives an important (Toeplitz)
groupoid
\[
\mathcal{T}_{m}:={\mathbb{Z}}^{m}\times\overline{{\mathbb{Z}}}^{m}%
|_{\overline{{\mathbb{Z}}}_{\geq}^{m}}=\{(j,k)\in{\mathbb{Z}}^{m}%
\times\overline{{\mathbb{Z}}}_{\geq}^{m}\;|\;j+k\in\overline{{\mathbb{Z}}%
}_{\geq}^{m}\}
\]
where $\overline{{\mathbb{Z}}}={\mathbb{Z}}\cup\{+\infty\}$ and $\overline
{{\mathbb{\ Z}}}_{\geq}:=\{0,1,2,3,...\}\cup\{+\infty\}$.

Recall that the C*-algebra $C(SU(n+1)_{q})$ is generated by $u_{ij}$, $1\leq
i,j\leq n+1$, satisfying $u^{\ast}u=uu^{\ast}=I$ and some other relations
\cite{Wo:sun,So}. In particular, the C*-algebra $C(SU(2)_{q})$ is generated by
$u_{ij}$, with $1\leq i,j\leq2$, satisfying $u_{22}=u_{11}^{\ast}$,
$u_{12}=-q^{-1}u_{21}^{\ast}$, and $u^{\ast}u=uu^{\ast}=1$
\cite{Wo:ts,VaSo:af}. An important irreducible (non-faithful) *-representation
$\pi_{0}$ of $C(SU(2)_{q})$, $q>1$, on $\ell^{2}({\mathbb{Z}}_{\geq})$ is
given by
\[
\pi_{0}(u)=\left(
\begin{array}
[c]{cc}%
\alpha & -q^{-1}\gamma\\
\gamma & \alpha^{\ast}%
\end{array}
\right)
\]
where $\alpha(e_{j})=(1-q^{-2j})^{1/2}e_{j-1}$ and $\gamma(e_{j})=q^{-j}e_{j}$
for $j\geq0$. Here $\pi_{0}$ is applied to $u=(u_{ij}) $ entrywise.

Irreducible one-dimensional *-representations of $C(SU(n+1)_{q})$ are defined
by $\tau_{t}(u_{ij})=\delta_{ij}t_{j}$ for $t\in\mathbb{T}^{n}$ (with
$t_{n+1}=t_{1}^{-1}t_{2}^{-1}...t_{n}^{-1}$). The $\mathbb{T}^{n}$-family
$\left\{  \tau_{t}\right\}  _{t\in\mathbb{T}^{n}}$ of one-dimensional
irreducible *-representations of $C\left(  SU(n+1)_{q}\right)  $ can be viewed
as a C*-algebra homomorphism $\tau_{n+1}=\tau:C\left(  SU(n+1)_{q}\right)
\rightarrow C\left(  \mathbb{T}^{n}\right)  \cong C^{\ast}\left(
\mathbb{Z}^{n}\right)  $. There are $n$ fundamental *-representations $\pi
_{i}=\pi_{0}\phi_{i}$ with $\phi_{i}:C(SU(n+1)_{q})\rightarrow C(SU(2)_{q})$
given by $\phi_{i}(u_{jk})=u_{j-i+1,k-i+1}$ if $\{j,k\}\subseteq\{i,i+1\}$ and
$\phi_{i}(u_{jk})=\delta_{jk}$ if otherwise.

The unique maximal element in the Weyl group of $SU\left(  n+1\right)  $ can
be expressed in the reduced form
\[
s_{1}s_{2}s_{1}s_{3}s_{2}s_{1}...s_{n}s_{n-1}...s_{2}s_{1}.
\]
So by Soibelman's classification of irreducible *-representations \cite{So},
$C(SU(n+1)_{q})$ can be embedded into
\[
C^{\ast}(\mathcal{G}^{n})\subseteq C^{\ast}({\mathbb{Z}}^{n})\otimes
\mathcal{B}(\ell^{2}({\ \mathbb{Z}}_{\geq}^{N}))
\]
by
\[
(\tau_{n+1}\otimes\pi_{121321...n(n-1)..21})\Delta^{N}%
\]
where $N=n(n+1)/2$,
\[
\pi_{i_{1}i_{2}...i_{m}}:=\pi_{i_{1}}\bigotimes\pi_{i_{2}}\bigotimes
...\bigotimes\pi_{i_{m}},
\]
$\Delta$ is the comultiplication on $C(SU(n+1)_{q})$ with $\Delta^{m}$ defined
recursively as $\Delta^{k}=\left(  \Delta\bigotimes\operatorname{id}\right)
\Delta^{k-1}$, and $\mathcal{G}^{n}$ is the groupoid ${\mathbb{Z}}^{n}%
\times{\mathbb{Z}}^{N}\times\overline{{\mathbb{Z}}}^{N}|_{\overline
{{\mathbb{Z}}}_{\geq}^{N}}$ with ${\mathbb{Z}}^{n}$ acting trivially on
$\overline{\mathbb{Z}}^{N}$, and ${\mathbb{Z}}^{N}$ acts by translation on
$\overline{{\mathbb{Z}}}^{N}$.

\section{Quantum spheres $\mathbb{S}_{q}^{2n+1}$}

Recall that the C*-algebra of the quantum sphere $S_{q}^{2n+1}=SU(n)_{q}%
\backslash SU(n+1)_{q}$ (as a homogeneous quantum space) is defined by
\[
C(S_{q}^{2n+1})=\{f\in C(SU(n+1)_{q}):(\Phi\otimes id)(\Delta f)=1\otimes f\}
\]
where $\Phi:C(SU(n+1)_{q})\rightarrow C(SU(n)_{q})$ is the quotient map
defined by $\Phi(u_{ij})=u_{ij}$, $\Phi(u_{n+1,k})=\Phi(u_{k,n+1})=0$, and
$\Phi(u_{n+1,n+1})=1$, for $1\leq i,j,k\leq n$. From
\[
u_{n+1,m}=(h_{n}\otimes1)(1\otimes u_{n+1,m})=(h_{n}\otimes1)(\Phi
\otimes1)(\Delta(u_{n+1,m})),
\]
we have
\[
u_{n+1,m}\in(h_{n}\otimes1)(\Phi\otimes1)\Delta(C(SU(n+1)_{q}))=C(S_{q}%
^{2n+1})
\]
by Nagy's result \cite{N:hm}, where $h_{n}$ is the Haar functional on
$C(SU(n)_{q}))$ \cite{Wo:cm,Sh:cqg} such that $h_{n}(1)=1$.

As discussed in \cite{Sh:cqg}, we can prove that $C(S_{q}^{2n+1})$ is
generated by $u_{n+1,m}$'s, by verifying that the monomials
\[
p^{i,j,k}=(u_{n+1,n}^{\ast})^{i_{n}}...(u_{n+1,1}^{\ast})^{i_{1}}y_{1}^{j_{1}%
}...y_{n}^{j_{n}}y_{n+1}^{j_{n+1}}(u_{n+1,1})^{k_{1}}...(u_{n+1,n})^{k_{n}}%
\]
are linearly independent, where $y_{m}=u_{n+1,m}u_{n+1,m}^{\ast}$ and
$i_{m},j_{m},k_{m}\geq0$ with $i_{m}k_{m}=0$ (this last condition was
accidentally missing in \cite{Sh:cqg}) for $1\leq m\leq n$, $y_{n+1}^{j_{n+1}%
}=(u_{n+1,n+1})^{j_{n+1}}$ if $j_{n+1}\geq0$, and $y_{n+1}^{j_{n+1}%
}=(u_{n+1,n+1}^{\ast})^{-j_{n+1}}$ if $j_{n+1}\leq0$.

So we have \textbf{\ }
\[
C(S_{q}^{2n+1})=C^{\ast}(\{u_{n+1,m}|\,1\leq m\leq n+1\})\subset C\left(
SU(n+1)_{q}\right)
\]
(used as the definition of $S_{q}^{2n+1}$ in \cite{VaSo:af}) which can be
embeded into the groupoid C*-algebra $C^{\ast}(\mathcal{F}^{n})$, where
$\mathcal{F}^{n}$ is the augmented $n$-dimensional Toeplitz groupoid
$\mathbb{Z}\times\mathcal{T}_{n}={\mathbb{Z}}\times({\mathbb{Z}}^{n}%
\times\overline{{\mathbb{Z}}}^{n}|_{\overline{{\mathbb{Z}}}_{\geq}^{n}})$,
through the faithful representation
\[
(\tau_{n+1}^{\prime}\otimes\pi_{n(n-1)..21})\Delta^{n}%
\]
where $\tau_{k}^{\prime}\left(  u_{ij}\right)  =\delta_{ki}\delta
_{ij}\operatorname{id}_{\mathbb{T}}\in C\left(  \mathbb{T}\right)  \cong
C^{\ast}\left(  \mathbb{Z}\right)  $ for $1\leq k\leq n+1$.

It turns out that $C(S_{q}^{2n+1})$ is actually isomorphic to the C*-algebra
of a subquotient groupoid of $\mathcal{F}^{n}$ \cite{Sh:qs}. In fact, let%

\[
\widetilde{\frak{F}_{n}}:=\{\left(  z,x,w\right)  \in\mathcal{F}^{n}%
|\;w_{i}=\infty\Longrightarrow
\]%
\[
x_{i}=-z-x_{1}-x_{2}-...-x_{i-1}\text{ and }x_{i+1}=...=x_{n}=0\}
\]
be a subgroupoid of $\mathcal{F}^{n}$. Define $\frak{F}_{n}:=\widetilde
{\frak{F}_{n}}/\sim$ where $\sim$ is the equivalence relation generated by
\[
(z,x,w)\sim(z,x,w_{1},...,w_{i}=\infty,\infty,...,\infty)
\]
for all $(z,x,w)$ with $w_{i}=\infty$ for an $1\leq i\leq n$. Then we have
$C(S_{q}^{2n+1})\simeq C^{\ast}(\frak{F}_{n})$.

Since $SU\left(  n\right)  $ is a Poisson Lie subgroup of $SU\left(
n+1\right)  $, the multiplicative Poisson structure on $SU\left(  n+1\right)
$ induces a natural covariant Poisson structure on the homogeneous space
$\mathbb{S}^{2n+1}=SU(n)\backslash SU(n+1)$, i.e. the $SU\left(  n+1\right)
$-action
\[
\mathbb{S}^{2n+1}\times SU(n+1)\rightarrow\mathbb{S}^{2n+1}%
\]
on $\mathbb{S}^{2n+1}$ is a Poisson map \cite{LuWe:plg}. Endowed with this
covariant Poisson structure, the sphere $\mathbb{S}^{2n+1}$ decomposes into
symplectic leaves \cite{We} which are the symplectic leaves in the canonically
embedded Poisson $\mathbb{S}^{2n-1}$ plus a circle family of leaves
symplectically isomorphic to the complex space $\mathbb{C}^{n}$. Analyzing the
above groupoid structure, we can easily see that this Poisson foliation
structure is reflected on the quantum algebras, namely, there is a short exact
sequence of C*-algebras
\[
0\rightarrow C\left(  \mathbb{T}\right)  \otimes\mathcal{K}\left(  \ell
^{2}\left(  \mathbb{Z}_{\geq}^{n}\right)  \right)  \rightarrow C(S_{q}%
^{2n+1})\rightarrow C(S_{q}^{2n-1})\rightarrow0.
\]
In fact, $\frak{F}_{n-1}\ $ can be identified with
\[
\left\{  \left[  (z,x,w)\right]  :(z,x,w)\in\widetilde{\frak{F}_{n}}\text{ and
}w_{n}=\infty\right\}  ,
\]
i.e. $\frak{F}_{n}$ restricted to the invariant closed subset $X$ consisting
of $\left[  w\right]  \in Z$ with $w_{n}=\infty$ in the unit space $Z$ of
$\frak{F}_{n}$, through the groupoid monomorphism sending $\left(
z,x^{\prime},w^{\prime}\right)  \in\widetilde{\frak{F}_{n-1}}$ to
\[
\left(  z,x^{\prime},-z-x_{1}^{\prime}-x_{2}^{\prime}-...-x_{n-1}^{\prime
},w^{\prime},\infty\right)  \in\widetilde{\frak{F}_{n}},
\]
and $\frak{F}_{n}\backslash\frak{F}_{n-1}=\frak{F}_{n}|_{Z\backslash X}$ (the
restriction of $\frak{F}_{n}$ to the invariant open set $Z\backslash X$) is
isomorphic to the groupoid ${\mathbb{Z}}^{1}\times\left(  {\mathbb{Z}}%
^{n}\times{\mathbb{Z}}^{n}|_{\mathbb{Z}_{\geq}^{n}}\right)  $. It is well
known \cite{Re} that there is a short exact sequence of groupoid C*-algebras
\[
0\rightarrow C^{\ast}\left(  \frak{F}_{n}\backslash\frak{F}_{n-1}\right)
\rightarrow C^{\ast}\left(  \frak{F}_{n}\right)  \rightarrow C^{\ast}\left(
\frak{F}_{n-1}\right)  \rightarrow0
\]
which leads to the above short exact sequence for the quantum spheres since
\[
C^{\ast}\left(  \frak{F}_{n}\backslash\frak{F}_{n-1}\right)  =C^{\ast}\left(
{\mathbb{Z}}^{1}\times\left(  {\mathbb{Z}}^{n}\times{\mathbb{Z}}%
^{n}|_{\mathbb{Z}_{\geq}^{n}}\right)  \right)  \cong C\left(  \mathbb{T}%
\right)  \otimes\mathcal{K}\left(  \ell^{2}\left(  \mathbb{Z}_{\geq}%
^{n}\right)  \right)  .
\]

\section{Standard quantum projective spaces $\mathbb{C}P_{q}^{n}$}

The C*-algebra of the standard quantum complex projective space%

\[
\mathbb{C}P_{q}^{n}:=U(n)_{q}\backslash SU(n+1)_{q},
\]
corresponding to the Poisson Lie subgroup $U\left(  n\right)  $ of $SU\left(
n+1\right)  $, can be identified with
\[
C(\mathbb{C}P_{q}^{n})=C^{\ast}(\{u_{n+1,i}^{\ast}u_{n+1,j}|\,1\leq i,j\leq
n+1\})\subset C\left(  \mathbb{S}_{q}^{2n+1}\right)  .
\]
In fact, as a quantum homogeneous space, $\mathbb{C}P_{q}^{n}$ is defined by
\[
C(\mathbb{C}P_{q}^{n})=\{f\in C(SU(n+1)_{q}):(\Phi^{\prime}\otimes id)(\Delta
f)=1\otimes f\}
\]
where $\Phi^{\prime}:C(SU(n+1)_{q})\rightarrow C(U(n)_{q})$ is the quotient
map defined by $\Phi(u_{ij})=u_{ij}$, $\Phi(u_{n+1,k})=\Phi(u_{k,n+1})=0$,
and
\[
\Phi(u_{n+1,n+1})=t_{n+1}=\left(  t_{1}...t_{n}\right)  ^{-1}\in C\left(
U\left(  n\right)  _{q}\right)  \subset C^{\ast}\left(  \mathbb{Z}^{n}%
\times\mathcal{T}_{n\left(  n-1\right)  /2}\right)
\]
for $1\leq i,j,k\leq n$. By Nagy's result \cite{N:hm}, we have
\[
(h_{n}^{\prime}\otimes1)(\Phi^{\prime}\otimes1)\Delta(C(SU(n+1)_{q}%
))=C(\mathbb{C}P_{q}^{n})
\]
where $h_{n}^{\prime}$ is the Haar functional on $C(U(n)_{q}))$ \cite{Wo:cm}
such that $h_{n}^{\prime}(1)=1$. From
\[
(h_{n}^{\prime}\otimes1)(\Phi^{\prime}\otimes1)(\Delta(u_{n+1,i}^{\ast
}u_{n+1,j}))=(h_{n}^{\prime}\otimes1)(\overline{t_{n+1}}t_{n+1}\otimes
u_{n+1,i}^{\ast}u_{n+1,j})
\]%
\[
=(h_{n}^{\prime}\otimes1)(1\otimes u_{n+1,i}^{\ast}u_{n+1,j})=u_{n+1,i}^{\ast
}u_{n+1,j},
\]
we get $u_{n+1,i}^{\ast}u_{n+1,j}\in C(\mathbb{C}P_{q}^{n})$ for $1\leq
i,j\leq n+1$. On the other hand, similar to the discussion in \cite{Sh:cqg},
we can prove that $C(\mathbb{C}P_{q}^{n})$ is generated by $u_{n+1,i}^{\ast
}u_{n+1,j}$, by verifying that the monomials
\[
P^{r,i,j,m}=z_{11}^{r_{1}}...z_{nn}^{r_{n}}z_{i_{1},j_{1}}...z_{i_{m},j_{m}}%
\]
(i) are linearly independent in the quantum case (i.e. when $q>1$) and (ii)
linearly span the *-algebra generated by $\overline{u_{n+1,i}}u_{n+1,j}$, or
the *-subalgebra of $\mathbb{T}$-invariant polynomials in $\overline
{u_{n+1,i}}$'s and $u_{n+1,j}$'s, in the classical case (i.e. when $q=1$),
where $m,r_{1},...,r_{n}\geq0$, $z_{i,j}=u_{n+1,i}^{\ast}u_{n+1,j}$, and the
two sets of indices $i_{1}\geq i_{2}\geq...\geq i_{m}$ and $j_{1}\leq
j_{2}\leq...\leq j_{m}$ are disjoint subsets of $\left\{  1,2,...,n+1\right\}
$. Note that if we allow all $u_{n+1,i}^{\ast}$'s and $u_{n+1,j}$'s in
$P^{i,j,m}$ to commute, then we can rewrite it, by simply permuting such
factors, in the form
\[
p^{i,j,k}=(u_{n+1,n}^{\ast})^{i_{n}}...(u_{n+1,1}^{\ast})^{i_{1}}y_{1}^{j_{1}%
}...y_{n}^{j_{n}}y_{n+1}^{j_{n+1}}(u_{n+1,1})^{k_{1}}...(u_{n+1,n})^{k_{n}}%
\]
where $y_{m}=u_{n+1,m}u_{n+1,m}^{\ast}$ and $i_{m},j_{m},k_{m}\geq0$ with
$i_{m}k_{m}=0$ for $1\leq m\leq n$, $y_{n+1}^{j_{n+1}}=(u_{n+1,n+1})^{j_{n+1}}
$ if $j_{n+1}\geq0$, and $y_{n+1}^{j_{n+1}}=(u_{n+1,n+1}^{\ast})^{-j_{n+1}} $
if $j_{n+1}\leq0$. So condition (ii) can be verified easily by integrating
$\mathbb{T}$-action on polynomials $p^{i,j,k}$ which form a linear basis of
polynomials in $\overline{u_{n+1,i}}$'s and $u_{n+1,j}$'s on $\mathbb{S}%
^{2n+1}$. On the other hand, condition (i) can be verified by considering the
weight functions and the orders of the weighted shifts of
\[
(\tau_{n+1}^{\prime}\otimes\pi_{n(n-1)..21})\Delta^{n}(P^{r,i,j,m})
\]
which can be determined by
\[
(\tau_{n+1}^{\prime}\otimes\pi_{n(n-1)..21})\Delta^{n}(u_{n+1,i}%
)=\operatorname{id}_{\mathbb{T}}\otimes\underset{n+1-i}{\underbrace
{\gamma\otimes...\otimes\gamma}}\otimes\underset{i-1}{\underbrace{\alpha
^{\ast}\otimes1\otimes...\otimes1}}.
\]

Under the faithful *-representation $\left(  \tau_{n+1}^{\prime}\otimes
\pi_{n\left(  n-1\right)  ...21}\right)  \Delta^{n}$, we have
\[
(\tau_{n+1}^{\prime}\otimes\pi_{n(n-1)..21})\Delta^{n}\left(  u_{n+1,i}^{\ast
}u_{n+1,j}\right)
\]%
\[
=\left[  (\tau_{n+1}^{\prime}\otimes\pi_{n(n-1)..21}\Delta^{n-1})\left(
\Delta u_{n+1,i}\right)  \right]  ^{\ast}(\tau_{n+1}^{\prime}\otimes
\pi_{n(n-1)..21}\Delta^{n-1})\left(  \Delta u_{n+1,j}\right)
\]%
\[
=(\operatorname{id}_{\mathbb{T}}\otimes\pi_{n(n-1)..21}\Delta^{n-1}%
u_{n+1,i})^{\ast}(\operatorname{id}_{\mathbb{T}}\otimes\pi_{n(n-1)..21}%
\Delta^{n-1}u_{n+1,j})
\]%
\[
=1\otimes\left(  \pi_{n(n-1)..21}\Delta^{n-1}u_{n+1,i}\right)  ^{\ast}%
(\pi_{n(n-1)..21}\Delta^{n-1}u_{n+1,j})
\]%
\[
=1\otimes\pi_{n(n-1)..21}\Delta^{n-1}(u_{n+1,i}^{\ast}u_{n+1,j}).
\]
So it is clear that the *-representation
\[
\pi_{n\left(  n-1\right)  ...21}\Delta^{n-1}=\left(  \pi_{n}\otimes\pi
_{n-1}\otimes...\otimes\pi_{1}\right)  \Delta^{n-1}%
\]
gives an embedding of $C(\mathbb{C}P_{q}^{n})$ into the groupoid C*-algebra
$C^{\ast}(\mathcal{T}_{n})$, where
\[
\mathcal{T}_{n}={\mathbb{Z}}^{n}\times\overline{{\mathbb{Z}}}^{n}%
|_{\overline{{\mathbb{Z}}}_{\geq}^{n}}.
\]

Applying a similar analysis as used in \cite{Sh:qs}, we can get the following
results. Let%

\[
\widetilde{\frak{T}_{n}}:=\{(j,k)\in\mathcal{T}_{n}|\;k_{i}=\infty
\Longrightarrow
\]%
\[
j_{i}=-j_{1}-j_{2}-...-j_{i-1}\text{ and }j_{i+1}=...=j_{n}=0\}
\]
be a subgroupoid of $\mathcal{T}_{n}$. Define a subquotient groupoid
$\frak{T}_{n}:=\widetilde{\frak{T}_{n}}/\sim$ where $\sim$ is the equivalence
relation generated by
\[
(j,k)\sim(j,k_{1},...,k_{i}=\infty,\infty,...,\infty)
\]
for all $(j,k)$ with $k_{i}=\infty$ for an $1\leq i\leq n$.

\begin{theorem}
$C(\mathbb{C}P_{q}^{n})\simeq C^{\ast}(\frak{T}_{n})$ and hence is independent
of $q$.
\end{theorem}

\begin{corollary}
There is a short exact sequence
\[
0\rightarrow\mathcal{K}\rightarrow C(\mathbb{C}P_{q}^{k})\rightarrow
C(\mathbb{C}P_{q}^{k-1})\rightarrow0
\]
for $k\geq1$ with $C(\mathbb{C}P_{q}^{0})\simeq\mathbb{C}$.
\end{corollary}

\begin{corollary}
The C*-algebra $C(\mathbb{C}P_{q}^{n})$ has the following composition
sequence,
\[
C(\mathbb{C}P_{q}^{n})=\mathcal{I}_{0}\supseteq\mathcal{I}_{1}\supseteq
...\supseteq\mathcal{I}_{n}\supseteq\mathcal{I}_{n+1}:=0,
\]
with
\[
\mathcal{I}_{k}/\mathcal{I}_{k+1}\simeq\mathcal{K}(\ell^{2}(\mathbb{Z}^{k}))
\]
for $k>0$ and $\mathcal{I}_{0}/\mathcal{I}_{1}\simeq\mathbb{C}$.
\end{corollary}

The above short exact sequences reflect faithfully the underlying singular
foliation defined by the canonical $SU\left(  n+1\right)  $-covariant Poisson
structure on $\mathbb{C}P^{n}=U\left(  n\right)  \backslash SU\left(
n+1\right)  $.

\section{Nonstandard quantum projective spaces $\mathbb{C}P_{q,c}^{n}$}

Recently Korogodsky and Vaksman studied nonstandard quantum projective spaces
$\mathbb{C}P_{q,c}^{n}$ \cite{KoVa}. Here we use the approach and result of
Dijkhuizen and Noumi \cite{DiNo}. They define $\mathbb{C}P_{q,c}^{n}$ by
\[
C\left(  \mathbb{C}P_{q,c}^{n}\right)  :=\left\{  v\in C\left(  SU(n+1)_{q}%
\right)  :v\cdot\frak{k}^{c}=0\right\}
\]
for some coideal $\frak{k}^{c}$ of the quantized universal enveloping algebra
$\mathcal{U}_{q}\left(  \frak{su}\left(  n+1\right)  \right)  $, $c\in\left(
0,\infty\right)  $, and show that
\[
C(\mathbb{C}P_{q,c}^{n})\cong C^{\ast}(\{x_{i}^{\ast}x_{j}|\,1\leq i,j\leq
n+1\})\subset C\left(  SU(n+1)_{q}\right)
\]
where
\[
x_{i}=\sqrt{c}u_{1,i}+u_{n+1,i}.
\]

We remark that when $c=0$, we have $C(\mathbb{C}P_{q,c}^{n})\cong
C(\mathbb{C}P_{q}^{n})$ the standard case, but for $c\in\left(  0,\infty
\right)  $, $\mathbb{C}P_{qc}^{n}$ does not come from $SU(n+1)_{q}$ modulo a
quantum subgroup.

To embed $C(\mathbb{C}P_{q,c}^{n})$ into a groupoid C*-algebra, we first note
that the unique maximal element in the Weyl group of $SU\left(  n+1\right)  $
can also be expressed in the reduced form
\[
s_{2}s_{3}s_{2}s_{4}s_{3}s_{2}...s_{n}s_{n-1}...s_{2}s_{1}s_{2}...s_{n-1}%
s_{n}.
\]
So by Soibelman's classification of irreducible *-representations \cite{So},
\[
(\tau_{n+1}\otimes\pi_{232432...n\left(  n-1\right)  ..2123...\left(
n-1\right)  n})\Delta^{N}%
\]
is a faithful *-representation of $C(SU(n+1)_{q})$.

Note that
\begin{align*}
&  (\tau_{n+1}\otimes\pi_{232432...n\left(  n-1\right)  ..2123...\left(
n-1\right)  n})\Delta^{N}\left(  u_{k,i}\right) \\
&  =t_{k}\otimes1\otimes...\otimes1\otimes\left(  \pi_{n\left(  n-1\right)
..2123...\left(  n-1\right)  n}\Delta^{2n-2}\left(  u_{k,i}\right)  \right)
\end{align*}
when $k=1$ or $k=n+1$. Thus $\left(  \tau_{n+1}\otimes\pi_{n\left(
n-1\right)  ..2123...\left(  n-1\right)  n}\right)  \Delta^{2n-1}$ is a
faithful *-representation of $C^{\ast}\left(  \left\{  u_{1,i},u_{n+1,i}%
\right\}  _{i=1}^{n+1}\right)  \supset C(\mathbb{C}P_{qc}^{n}) $. Furthermore,
since
\[
\left(  \tau_{n+1}\otimes\pi_{n\left(  n-1\right)  ..2123...\left(
n-1\right)  n}\right)  \Delta^{2n-1}\left(  x_{i}\right)  =t_{1}\otimes
\sqrt{c}\pi_{n\left(  n-1\right)  ..2123...\left(  n-1\right)  n}\Delta
^{2n-2}\left(  u_{1,i}\right)
\]%
\[
+\left(  t_{1}t_{2}...t_{n}\right)  ^{-1}\otimes\pi_{n\left(  n-1\right)
..2123...\left(  n-1\right)  n}\Delta^{2n-2}\left(  u_{n+1,i}\right)
\]%
\[
\in\left\{  t_{1}^{l}\left(  t_{2}...t_{n}\right)  ^{m}:l-2m=1\right\}
\otimes C^{\ast}\left(  {\mathbb{Z}}^{2n-1}\times\overline{{\mathbb{Z}}%
}^{2n-1}|_{\overline{{\mathbb{Z}}}_{\geq}^{2n-1}}\right)  ,
\]
we have
\[
\left(  \tau_{n+1}\otimes\pi_{n\left(  n-1\right)  ..2123...\left(
n-1\right)  n}\right)  \Delta^{2n-1}\left(  x_{i}^{\ast}x_{j}\right)
\]%
\[
\in\left\{  t_{1}^{l}\left(  t_{2}...t_{n}\right)  ^{m}:l-2m=0\right\}
\otimes C^{\ast}\left(  {\mathbb{Z}}^{2n-1}\times\overline{{\mathbb{Z}}%
}^{2n-1}|_{\overline{{\mathbb{Z}}}_{\geq}^{2n-1}}\right)
\]%
\[
=\left\{  \left(  t_{1}^{2}t_{2}...t_{n}\right)  ^{m}:m\in\mathbb{Z}\right\}
\otimes C^{\ast}\left(  {\mathbb{Z}}^{2n-1}\times\overline{{\mathbb{Z}}%
}^{2n-1}|_{\overline{{\mathbb{Z}}}_{\geq}^{2n-1}}\right)  .
\]
Thus $\left(  \tau_{n+1}\otimes\pi_{n\left(  n-1\right)  ..2123...\left(
n-1\right)  n}\right)  \Delta^{2n-1}$ embeds $C(\mathbb{C}P_{qc}^{n})$ into
\[
C^{\ast}\left(  \mathbb{Z}\right)  \otimes C^{\ast}\left(  {\mathbb{Z}}%
^{2n-1}\times\overline{{\mathbb{Z}}}^{2n-1}|_{\overline{{\mathbb{Z}}}_{\geq
}^{2n-1}}\right)
\]%
\[
\cong C^{\ast}\left(  \left\{  t_{1}^{2}t_{2}...t_{n}\right\}  \right)
\otimes C^{\ast}\left(  {\mathbb{Z}}^{2n-1}\times\overline{{\mathbb{Z}}%
}^{2n-1}|_{\overline{{\mathbb{Z}}}_{\geq}^{2n-1}}\right)
\]%
\[
\subset C\left(  \mathbb{T}^{n}\right)  \otimes C^{\ast}\left(  {\mathbb{Z}%
}^{2n-1}\times\overline{{\mathbb{Z}}}^{2n-1}|_{\overline{{\mathbb{Z}}}_{\geq
}^{2n-1}}\right)  .
\]

\begin{theorem}
For $0<c<\infty$, $C(\mathbb{C}P_{qc}^{n})$ can be embedded into the groupoid
C*-algebra $C^{\ast}\left(  \frak{G}_{n}\right)  $ where
\[
\frak{G}_{n}:={\mathbb{Z}}\times\left(  {\mathbb{Z}}^{2n-1}\times
\overline{{\mathbb{Z}}}^{2n-1}|_{\overline{{\mathbb{Z}}}_{\geq}^{2n-1}%
}\right)  .
\]
\end{theorem}

When $n=1$, we have $\mathbb{C}P^{1}=\mathbb{S}^{2}$ in the classical case,
while in the quantum case, $\mathbb{C}P_{qc}^{1}$ is indeed a Podle\`{s}'
quantum sphere $\mathbb{S}_{\mu c}^{2}$ \cite{Po} (with $\mu=q^{-1}$) for some
(different) $c$. In \cite{LuWe:ps} Lu and Weinstein classified all $SU\left(
2\right)  $-covariant Poisson structures on $\mathbb{S}^{2}$ by a real
parameter and showed that each `nonstandard' $SU\left(  2\right)  $-covariant
Poisson sphere (corresponding to a suitable parameter $c\in\left(
0,\infty\right)  $) contains the trivial Poisson $1$-sphere $\mathbb{S}^{1}$
(consisting of a circle family of $0$-dimensional symplectic leaves) and
exactly two $2$-dimensional symplectic leaves. This geometric structure is
again reflected faithfully in the algebraic structure of the algebra $C\left(
\mathbb{S}_{\mu c}^{2}\right)  $ of the nonstandard quantum spheres
$\mathbb{S}_{\mu c}^{2}$ as follows.

\begin{theorem}
For $0<c<\infty$, $C(\mathbb{S}_{\mu c}^{2})\cong C^{\ast}(\frak{G}^{\prime}%
)$, where
\[
\frak{G}^{\prime}=\{(j,j,k_{1},k_{2})\;|\;\text{ }k_{1}=\infty\text{ or }%
k_{2}=\infty\}
\]
is a subgroupoid of the groupoid ${\mathbb{Z}}^{2}\times\overline{{\mathbb{Z}%
}}^{2}|_{\overline{{\mathbb{Z}}}_{\geq}^{2}}$, and there is a short exact
sequence
\[
0\rightarrow\mathcal{K}\oplus\mathcal{K}\rightarrow C(\mathbb{S}_{\mu c}%
^{2})\rightarrow C(\mathbb{S}^{1})\rightarrow0.
\]
\end{theorem}

On higher dimensional projective spaces $\mathbb{C}P^{n}$, we also have a
one-parameter family of nonstandard $SU\left(  n+1\right)  $-covariant Poisson
structures suitably parametrized by $c\in\left(  0,\infty\right)  $, and
similar to the case of $n=1$, each such nonstandard Poisson $\mathbb{C}P^{n}$
contains an embedded copy of a standard Poisson $\mathbb{S}^{2n-1}$
\cite{Sh:cp}. On the quantum level, one would then expect that $C(\mathbb{C}%
P_{qc}^{n})$ should have $C\left(  \mathbb{S}_{q}^{2n-1}\right)  $ as a
quotient to reflect this geometric fact. Using the above groupoid description
of $C(\mathbb{C}P_{qc}^{n})$, we can indeed show that this is the case.

\begin{theorem}
For $0<c<\infty$, there is a short exact sequence
\[
0\rightarrow\mathcal{I}\rightarrow C(\mathbb{C}P_{qc}^{n})\rightarrow C\left(
\mathbb{S}_{q}^{2n-1}\right)  \rightarrow0
\]
for some ideal $\mathcal{I}$.
\end{theorem}

\begin{proof}
First we note that since $\left(  \tau_{n+1}\otimes\pi_{n\left(  n-1\right)
..212...\left(  n-1\right)  n}\right)  \Delta^{2n-1}$ is faithful on $C\left(
\mathbb{S}_{q}^{2n+1}\right)  $,
\[
\left(  \tau_{n+1}\otimes\pi_{n\left(  n-1\right)  ..212...\left(  n-1\right)
n}\right)  \Delta^{2n-1}\left(  u_{n+1,i}\right)  =t_{n+1}\otimes\pi_{n\left(
n-1\right)  ..2123...\left(  n-1\right)  n}\Delta^{2n-2}u_{n+1,i}%
\]%
\[
=\left(  1-\delta_{n+1,i}\right)  t_{n+1}\otimes\gamma^{\otimes n-i}%
\otimes\underset{i-1}{\underbrace{\alpha^{\ast}\otimes1\otimes...\otimes1}%
}\otimes1^{\otimes i-1}\otimes\underset{\left(  n+1\right)  -i}{\underbrace
{\gamma\otimes1\otimes...\otimes1}}+\left(  1-\delta_{1,i}\right)  \left(
1-\delta_{2,i}\right)
\]%
\[
\left[  \sum_{k=n+1-i}^{n-2}t_{n+1}\otimes\gamma^{\otimes k}\otimes
\alpha^{\ast}\otimes1^{\otimes n-k-1}\otimes1^{\otimes n-k-2}\otimes
\alpha^{\ast}\otimes\left(  -q^{-1}\gamma\right)  ^{\otimes k-\left(
n+1\right)  +i}\otimes\underset{\left(  n+1\right)  -i}{\underbrace
{\alpha\otimes1\otimes...\otimes1}}\right]
\]%
\[
+\left(  1-\delta_{1,i}\right)  t_{n+1}\otimes\gamma^{\otimes n-1}%
\otimes\alpha^{\ast}\otimes\left(  -q^{-1}\gamma\right)  ^{\otimes i-2}%
\otimes\underset{\left(  n+1\right)  -i}{\underbrace{\alpha\otimes
1\otimes...\otimes1}}%
\]
with $1\leq i\leq n+1$ generate a C*-algebra isomorphic to $C\left(
\mathbb{S}_{q}^{2n+1}\right)  $.

Now let $\rho:C^{\ast}\left(  {\mathbb{Z}}\times\overline{{\mathbb{Z}}%
}|_{\overline{{\mathbb{Z}}}_{\geq}}\right)  \rightarrow\mathbb{C}$ be the
composition of the C*-homomorphism
\[
C^{\ast}\left(  {\mathbb{Z}}\times\overline{{\mathbb{Z}}}|_{\overline
{{\mathbb{Z}}}_{\geq}}\right)  \rightarrow C^{\ast}\left(  {\mathbb{Z}}%
\times\overline{{\mathbb{Z}}}|_{\left\{  \infty\right\}  }\right)  =C^{\ast
}\left(  {\mathbb{Z}}\right)  \cong C\left(  \mathbb{T}\right)
\]
induced by the restriction to the invariant closed subset $\left\{
\infty\right\}  $ of the unit space $\overline{{\mathbb{Z}}}_{\geq}$ and the
evaluation map $C\left(  \mathbb{T}\right)  \rightarrow\mathbb{C}$ at
$1\in\mathbb{T}$, and let
\[
\tilde{\rho}:=\underset{n-1}{\underbrace{1\otimes...\otimes1}}\otimes
\rho\otimes\rho\otimes\underset{n-2}{\underbrace{1\otimes...\otimes1}}%
:C^{\ast}\left(  \mathcal{T}_{2n-1}\right)  \rightarrow C^{\ast}\left(
\mathcal{T}_{2n-3}\right)
\]
be the homomorphism that `removes' the middle two tensor factors. Then by
direct calculation, we get
\[
\left(  1\otimes\tilde{\rho}\right)  \left(  \tau_{n+1}\otimes\pi_{n\left(
n-1\right)  ..212...\left(  n-1\right)  n}\right)  \Delta^{2n-1}\left(
\sqrt{c}u_{1,1}+u_{n+1,1}\right)
\]%
\[
=\sqrt{c}t_{1}\otimes1\otimes.....\otimes1
\]
while
\[
\left(  1\otimes\tilde{\rho}\right)  \left(  \tau_{n+1}\otimes\pi_{n\left(
n-1\right)  ..212...\left(  n-1\right)  n}\right)  \Delta^{2n-1}\left(
\sqrt{c}u_{1,i}+u_{n+1,i}\right)
\]%
\[
=t_{n+1}\otimes\pi_{n\left(  n-1\right)  ..323...\left(  n-1\right)  n}%
\Delta^{2n-4}\left(  u_{n+1,i}\right)
\]
for $i>1$. So
\[
\left(  1\otimes\tilde{\rho}\right)  \left(  \tau_{n+1}\otimes\pi_{n\left(
n-1\right)  ..212...\left(  n-1\right)  n}\right)  \Delta^{2n-1}\left(
x_{1}^{\ast}x_{1}\right)  =c,
\]
and for $i>1$,
\[
\left[  \left(  1\otimes\tilde{\rho}\right)  \left(  \tau_{n+1}\otimes
\pi_{n\left(  n-1\right)  ..212...\left(  n-1\right)  n}\right)  \Delta
^{2n-1}\left(  x_{i}^{\ast}x_{1}\right)  \right]  ^{\ast}%
\]%
\[
=\left(  1\otimes\tilde{\rho}\right)  \left(  \tau_{n+1}\otimes\pi_{n\left(
n-1\right)  ..212...\left(  n-1\right)  n}\right)  \Delta^{2n-1}\left(
x_{1}^{\ast}x_{i}\right)
\]%
\[
=\sqrt{c}\left(  t_{1}^{2}t_{2}...t_{n}\right)  ^{-1}\otimes\pi_{n\left(
n-1\right)  ..323...\left(  n-1\right)  n}\Delta^{2n-4}\left(  u_{n+1,i}%
\right)
\]%
\[
=\sqrt{c}\left(  \tau^{\prime\prime}\otimes\pi_{n\left(  n-1\right)
..323...\left(  n-1\right)  n}\right)  \Delta^{2n-3}\left(  u_{n+1,i}\right)
\]
while for $i,j>1$,
\[
\left(  1\otimes\tilde{\rho}\right)  \left(  \tau_{n+1}\otimes\pi_{n\left(
n-1\right)  ..212...\left(  n-1\right)  n}\right)  \Delta^{2n-1}\left(
x_{i}^{\ast}x_{j}\right)
\]%
\[
=1\otimes\pi_{n\left(  n-1\right)  ..323...\left(  n-1\right)  n}\Delta
^{2n-4}\left(  u_{n+1,i}^{\ast}u_{n+1,j}\right)
\]%
\[
=\left(  \tau^{\prime\prime}\otimes\pi_{n\left(  n-1\right)  ..323...\left(
n-1\right)  n}\right)  \Delta^{2n-3}\left(  u_{n+1,i}^{\ast}u_{n+1,j}\right)
,
\]
where $\tau^{\prime\prime}\left(  u_{ij}\right)  =\delta_{i,j}\delta
_{n+1,i}\left(  t_{1}^{2}t_{2}...t_{n}\right)  ^{-1}$. It is straight forward
to verify that
\[
\left(  \tau^{\prime\prime}\otimes\pi_{n\left(  n-1\right)  ..323...\left(
n-1\right)  n}\right)  \Delta^{2n-3}\left(  u_{n+1,i}\right)  ,
\]
with $1<i\leq n+1$, coincides with the operators
\[
\left(  \tau_{n}\otimes\pi_{\left(  n-1\right)  ..212...\left(  n-1\right)
}\right)  \Delta^{2n-3}\left(  u_{n,i-1}\right)
\]
whose formula is explicitly written above for the case of dimension $n+1$. So
\[
\left(  1\otimes\tilde{\rho}\right)  \left(  \tau_{n+1}\otimes\pi_{n\left(
n-1\right)  ..212...\left(  n-1\right)  n}\right)  \Delta^{2n-1}\left(
C(\mathbb{C}P_{qc}^{n})\right)
\]%
\[
\cong\left(  \tau_{n}\otimes\pi_{\left(  n-1\right)  ..212...\left(
n-1\right)  }\right)  \Delta^{2n-3}\left(  C(\mathbb{S}_{q}^{2n-1})\right)
\cong C(\mathbb{S}_{q}^{2n-1})
\]
and hence $C(\mathbb{S}_{q}^{2n-1})$ is a quotient of $C(\mathbb{C}P_{qc}%
^{n})$.
\end{proof}

\end{document}